\documentclass[11pt]{article}
\usepackage{enumerate}
\usepackage{amssymb,a4wide,latexsym,makeidx,epsfig,fleqn}
\usepackage{amsthm}
\usepackage{amsmath}
\usepackage{enumerate}
\newtheorem{theorem}{Theorem}[section]
\newtheorem{remark}[theorem]{Remark}

\newtheorem{lemma}[theorem]{Lemma}

\newtheorem{corollary}[theorem]{Corollary}

\begin{document}
\textwidth 150mm \textheight 225mm
\title{Hermitian-Randi\'{c} matrix and Hermitian-Randi\'{c} energy of mixed graphs
\thanks{ Supported by
the National Natural Science Foundation of China (No. 11171273)}}
\author{{Yong Lu$^1$, Ligong Wang$^{2,}$\footnote{Corresponding author.} and Qiannan Zhou$^3$}\\
{\small Department of Applied Mathematics, School of Science, Northwestern
Polytechnical University,}\\ {\small  Xi'an, Shaanxi 710072,
People's Republic
of China.}\\
{\small $^1$ E-mail: luyong.gougou@163.com}\\
{\small $^2$ E-mail: lgwangmath@163.com}\\
{\small $^3$ E-mail: qnzhoumath@163.com}\\}
\date{}
\maketitle
\begin{center}
\begin{minipage}{120mm}
\vskip 0.3cm
\begin{center}
{\small {\bf Abstract}}
\end{center}
{\small Let $M$ be a mixed graph and $H(M)$ be its Hermitian-adjacency matrix. If we add every edge and arc in $M$ a Randi\'{c} weight, then we can get a new weighted Hermitian-adjacency matrix. What are the properties of this new matrix? Motivated by this, we define the Hermitian-Randi\'{c} matrix $R_{H}(M)=(r_{h})_{kl}$ of a mixed graph $M$,
where $(r_{h})_{kl}=-(r_{h})_{lk}=\frac{\textbf{\textrm{i}}}{\sqrt{d_{k}d_{l}}}$ ($\textbf{\textrm{i}}=\sqrt{-1}$) if $(v_{k},v_{l})$ is an arc of $M$,  $(r_{h})_{kl}=(r_{h})_{lk}=\frac{1}{\sqrt{d_{k}d_{l}}}$ if $v_{k}v_{l}$ is an undirected edge of $M$, and $(r_{h})_{kl}=0$ otherwise.
 In this paper, firstly, we compute the characteristic polynomial of the Hermitian-Randi\'{c} matrix of a mixed graph. Furthermore, we give bounds to the Hermitian-Randi\'{c} energy of a general mixed graph. Finally, we give some results about the Hermitian-Randi\'{c} energy of mixed trees.

 }
\end{minipage}
\end{center}

\section{Introduction }
 In this paper, we only consider simple graphs without multiedges and loops. A graph $M$ is said to be \emph{mixed} if it is obtained from an undirected graph $M_{U}$ by orienting a subset of its edges. We call $M_{U}$ the \emph{underlying graph} of $M$. Clearly, a mixed graph concludes both possibilities of all edges oriented and all edges undirected as extreme cases.

Let $M$ be a mixed graph with vertex set $V(M)=\{v_{1},v_{2},\ldots,v_{n}\}$ and edge set $E(M)$. For $v_{i},v_{j}\in V(M)$, we denote an undirected edge joining two vertices $v_{i}$ and $v_{j}$ of $M$ by $v_{i}v_{j}$ (or $v_{i}\leftrightarrow v_{j}$). Denote a directed edge (or arc) from $v_{i}$ to $v_{j}$ by $(v_{i},v_{j})$ (or $v_{i}\rightarrow v_{j}$). In addition, let $E_{0}(M)$ denote the set of all undirected edges and $E_{1}(M)$ denote all the directed arcs set. Clearly, $E(M)$ is the union of $E_{0}(M)$ and $E_{1}(M)$. A mixed graph is called \emph{mixed tree} (or \emph{mixed bipartite graph}) if its underlying graph is a tree (or bipartite graph). In general, the order, size, number of components, degree of a vertex of $M$ are the same to those in $M_{U}$. We use Bondy and Murty \cite{BM} for terminologies and notations not defined here.

Let $G$ be a  simple graph with vertex set $\{v_{1},v_{2},\ldots,v_{n}\}$. The \emph{adjacency matrix} of a simple graph $G$  of order $n$ is defined as the $n\times n$ symmetric square matrix $A=A(G)=(a_{ij})$, where  $a_{ij}=1$  if $v_{i}v_{j}$ is an edge of $G$, otherwise $a_{ij}=0$.
 We denote by $d_{i}=d(v_{i})=d_{G}(v_{i})~(i=1,2,\ldots,n)$ the degree of vertex $v_{i}$. In addition, for a mixed graph $M$, if $v_{i}\in V(M)$, then we also denote by $d_{i}=d(v_{i})=d_{M_{U}}(v_{i})$.  The \emph{energy} of the graph $G$ (See the survey of Gutman, Li and Zhang \cite{ILZ} and the book of Li, Shi and Gutman \cite{LSI}) is defined as $\mathcal{E}_{A}(G)=\sum_{i=1}^{n}|\rho_{i}|$, where  $\rho_{1},\rho_{2},\ldots,\rho_{n}$ are all eigenvalues of $A(G)$.

 A convenient parameter of $G$ is the \emph{general Randi\'{c} index} $R_{\alpha}(G)$ defined as $R_{\alpha}(G)=\sum_{uv\in E(G)}(d_{u}d_{v})^{\alpha}$, where the summation is over all (unordered) edges $uv$ in $G$.
The molecular structure-descriptor, first proposed by Randi\'{c} \cite{Randic} in 1975,
 is defined as the sum of $\frac{1}{\sqrt{d_{u}d_{v}}}$ over all edges $uv$ of $G$ (with $\alpha=-\frac{1}{2}$). Nowadays,  $R=R(G)=\sum_{uv\in E(G)}\frac{1}{\sqrt{d_{u}d_{v}}}$ of $G$ is referred to as the \emph{Randi\'{c} index}. Countless chemical applications, the mathematical properties and mathematical chemistry of the Randi\'{c} index were reported in \cite{LXIG, LIS, RandicH}.

Gutman et al.\cite{IBS} pointed out that the Randi\'{c}-index-concept is purposeful to associate the graph $G$ a symmetric square matrix of order $n$, named \emph{Randi\'{c} matrix} $R(G)=(r_{ij})$, where
  $r_{ij}=\frac{1}{\sqrt{d_{i}d_{j}}}$ if $v_{i}v_{j}$ is an edge of $G$, otherwise $r_{ij}=0$. Let $D(G)$ be the diagonal matrix of vertex degrees of $G$. If $G$ has no isolated vertices, then $R(G)=D(G)^{-\frac{1}{2}}A(G)D(G)^{-\frac{1}{2}}$.

  The concept of the \emph{Randi\'{c} energy} of a graph $G$, denoted by $\mathcal{E}_{R}(G)$, was introduced in \cite{SAIA} as $\mathcal{E}_{R}(G)=\sum_{i=1}^{n}|\gamma_{i}|$, where $\gamma_{i}$ is the eigenvalues of $R(G)$, $i=1,2,\ldots,n$. Some basic properties of the Randi\'{c} index, Randi\'{c} matrix and Randi\'{c} energy were determined in the papers \cite{SDR, SDS, SAI, SAIA, GHLG, GHL, GLL, IBS, LGS, LXWJ, LXYYB, LXYYS, SYT}.

An oriented graph $G^{\sigma}$ is a digraph which assigns each edge of $G$ with a direction $\sigma$. The \emph{skew adjacency matrix} associated to $G^{\sigma}$ is the $n\times n$ matrix $S(G^{\sigma})=(s_{ij})$, where $s_{ij}=-s_{ji}=1$  if $(v_{i},v_{j})$ is an arc of $G^{\sigma}$, otherwise $s_{ij}=s_{ji}=0$. The \emph{skew energy} of $G^{\sigma}$, denoted by $\mathcal{E}_{S}(G^{\sigma})$, is defined as the sum of the norms of all the eigenvalues of $S(G^{\sigma})$. For more details about skew energy we can refer to \cite{CRST, LXLH}.

In 2016, Gu, Huang and Li \cite{GHL} defined the \emph{skew Randi\'{c} matrix}  $R_{s}(G^{\sigma})=((r_{s})_{ij})$ of an oriented graph $G^{\sigma}$ of order $n$, where $(r_{s})_{ij}=-(r_{s})_{ji}=\frac{1}{\sqrt{d_{i}d_{j}}}$  if $(v_{i},v_{j})$ is an arc of $G^{\sigma}$, otherwise $(r_{s})_{ij}=(r_{s})_{ji}=0$. Let $D(G)$ be the diagonal matrix of vertex degrees of $G$. If $G^{\sigma}$ has no isolated vertices, then $R_{s}(G^{\sigma})=D(G)^{-\frac{1}{2}}S(G^{\sigma})D(G)^{-\frac{1}{2}}$.

The \emph{Hermitian-adjacency matrix} of a mixed graph $M$ of order $n$ is the $n\times n$ matrix $H(M)=(h_{kl})$, where $h_{kl}=-h_{lk}=\textbf{\textrm{i}}$ ($\textbf{\textrm{i}}=\sqrt{-1}$) if $(v_{k},v_{l})$ is an arc of $M$, $h_{kl}=h_{lk}=1$ if $v_{k}v_{l}$ is an undirected edge of $M$, and $h_{kl}=0$ otherwise.
Obviously, $H(M)=H(M)^{\ast}:={\overline{H(M)}}^{T}$. Thus all its eigenvalues are real. This matrix  was introduced by Liu and Li in \cite{LJLX} and independently by Guo and Mohar in \cite{GM}.
The \emph{Hermitian energy}  of a mixed graph $M$ is defined as $\mathcal{E}_{H}(M)=\sum^{n}_{i=1}|\lambda_{i}|$, where $\lambda_{1},\lambda_{2},\ldots,\lambda_{n}$ are all eigenvalues of $H(M)$. Denote by $Sp_{H}(M)=(\lambda_{1},\lambda_{2},\ldots,\lambda_{n})$ the \emph{spectrum} of $H(M)$. For more details about Hermitian-adjacency matrix and Hermitian energy of mixed graphs, we can refer to \cite{CLZ, GM, LJLX, BMO, QHYG, YQ}.

From the above we can see that if we add every edge in a simple graph $G$ a Randi\'{c} weight, then we can get a Randi\'{c} matrix $R(G)$. If we add every arc in an oriented graph $G^{\sigma}$ a Randi\'{c} weight, then we can get a skew Randi\'{c} matrix $R_{s}(G^{\sigma})$.
Let $M$ be a mixed graph and $H(M)$ be its Hermitian-adjacency matrix, if we add every edge and arc in $M$ a Randi\'{c} weight, then we can get a new weighted Hermitian-adjacency matrix. What are the properties of this new matrix?  Motivated by this, we define the Hermitian-Randi\'{c} matrix of a mixed graph $M$.

Let $M$ be a mixed graph on the vertex set $\{v_{1},v_{2},\ldots,v_{n}\}$, then the  \emph{Hermitian-Randi\'{c} matrix} of $M$ is the $n\times n$ matrix $R_{H}(M)=((r_{h})_{kl})$, where
\begin{displaymath}
(r_{h})_{kl}=\left\{\
        \begin{array}{ll}
          \frac{1}{\sqrt{d_{k}d_{l}}},&\textrm{if}~v_{k}\leftrightarrow v_{l},\\
         \frac{\textbf{\textrm{i}}}{\sqrt{d_{k}d_{l}}},&\textrm{if}~v_{k}\rightarrow v_{l},\\
         \frac{-\textbf{\textrm{i}}}{\sqrt{d_{k}d_{l}}},&\textrm{if}~v_{l}\rightarrow v_{k},\\
         0,&\textrm{otherwise}.
        \end{array}
      \right.
\end{displaymath}

Let $D(M_{U})$ be the diagonal matrix of vertex degrees of $M_{U}$.
If $M$ has no isolated vertices, then $R_{H}(M)=D(M_{U})^{-\frac{1}{2}}H(M)D(M_{U})^{-\frac{1}{2}}$.
For a mixed graph $M$, let $R_{H}(M)$ be its Hermitian-Randi\'{c} matrix. It is obvious that $R_{H}(M)$ is Hermitian matrix, so all its eigenvalues $\mu_{1},\mu_{2},\ldots,\mu_{n}$ are real. The \emph{spectrum} of $R_{H}(M)$ is defined as $Sp_{R_{H}}(M)=(\mu_{1},\mu_{2},\ldots,\mu_{n})$.  The energy of $R_{H}(M)$, denoted by $\mathcal{E}_{R_{H}}(M)$, called \emph{Hermitian-Randi\'{c} energy}, which is defined as the sum of the absolute values of its eigenvalues of $R_{H}(M)$, that is $\mathcal{E}_{R_{H}}(M)=\sum^{n}_{i=1}|\mu_{i}|$.

\section{Hermitian-Randi\'{c} characteristic polynomial of  a mixed graph}

  In this section, we will give the characteristic polynomial of Hermitian-Randi\'{c} matrix of a mixed graph $M$, i.e., the $R_{H}$-characteristic polynomial of $M$. At first, we will introduce some basic definitions.

The \emph{value} of a mixed walk $W=v_{1}v_{2}\cdots v_{l}$ is $r_{h}(W)=(r_{h})_{12}(r_{h})_{23}(r_{h})_{(l-1)l}$. A mixed walk $W$ is \emph{positive} (or \emph{negative}) if $r_{h}(W)=\frac{1}{\sqrt{d_{1}d_{l}}d_{2}d_{3}\cdots d_{(l-1)}}$ (or $r_{h}(W)=-\frac{1}{\sqrt{d_{1}d_{l}}d_{2}d_{3}\cdots d_{(l-1)}}$). Note that for one direction the value of a mixed walk or a mixed cycle is $\alpha$, then for the reversed direction its value is $\overline{\alpha}$. Thus, if the value of a mixed cycle $C$ is $\prod_{v_{j}\in V(C)}\frac{1}{d(v_{j})}$ (resp. $-\prod_{v_{j}\in V(C)}\frac{1}{d(v_{j})}$) in a direction, then its value is $\prod_{v_{j}\in V(C)}\frac{1}{d(v_{j})}$ (resp. $-\prod_{v_{j}\in V(C)}\frac{1}{d(v_{j})}$) for the reversed direction. In these situations, we just term this mixed cycle as a positive (resp. negative) mixed cycle without mentioning any direction.

If each mixed cycle is positive (resp. negative) in a  mixed graph $M$, then $M$ is \emph{positive} (resp. \emph{negative}). A mixed graph $M$ is called an \emph{elementary graph}  if every component of $M$ is an edge, an arc or a mixed cycle, where every edge-component in $M$ is defined to be positive.  A  \emph{real spanning elementary subgraph} of a mixed graph $M$ is an elementary subgraph such that it contains all vertices of $M$ and all its mixed cycles are real.

Now we will give two results which are similar to those of \cite{GX, HL, LJLX}.

Let $M$ be a mixed graph of order $n$ with its Hermitian-Randi\'{c} matrix $R_{H}(M)$. Denote the \emph{$R_{H}$-characteristic polynomial} of $R_{H}(M)$ of $M$ by
$$P_{R_{H}}(M,x)=\mathrm{det}(xI-R_{H}(M))=x^{n}+a_{1}x^{n-1}+a_{2}x^{n-2}+\cdots+a_{n}.$$

\noindent\begin{theorem}\label{th:2.1}
Let $R_{H}(M)$ be the Hermitian-Randi\'{c} matrix of a mixed graph $M$ of order $n$. Then
$$\mathrm{det}R_{H}(M)=\sum_{M'}(-1)^{r(M')+l(M')}2^{s(M')}W(M'),$$
where the summation is over all real spanning elementary subgraphs $M'$ of $M$, $r(M')=n-c(M')$, $c(M')$ denotes the number of components of $M'$, $l(M')$ denotes the number of negative mixed cycles of $M'$, $s(M')$ denotes the number of mixed cycles with length $\geq3$ in $M'$, $W(M')=\prod_{v_{i}\in V(M')}\frac{1}{d_{M_{U}}(v_{i})}$.

\end{theorem}

\noindent {\bf Proof.}
   Let $M$ be a mixed graph of order $n$ with vertex set $\{v_{1},v_{2},\ldots, v_{n}\}$. Then
$$\mathrm{det}R_{H}(M)=\sum_{\pi\in S_{n}}sgn(\pi)(r_{h})_{1\pi(1)}(r_{h})_{2\pi(2)}\cdots (r_{h})_{n\pi(n)},$$
where $S_{n}$ is the set of all permutations on $\{1,2,\ldots,n\}$.

Consider a term $sgn(\pi)(r_{h})_{1\pi(1)}(r_{h})_{2\pi(2)}\cdots (r_{h})_{n\pi(n)}$ in the expansion of $\mathrm{det}R_{H}(M)$. If $v_{k}v_{\pi(k)}$ is not an edge or arc of $M$, then $(r_{h})_{k\pi(k)}=0$; that is, this term vanishes. Thus, if the term corresponding to a permutation $\pi$ is non-zero, then $\pi$ is fixed-point-free and can be expressed uniquely as the composition of disjoint cycles of length at least 2. Consequently, each non-vanishing term in the expansion of $\mathrm{det}R_{H}(M)$ gives rise to an elementary mixed graph $M'$ of $M$, with $V(M')=V(M)$. That is, $M'$ is a spanning elementary subgraph of $M$ of order $n$.

A spanning elementary subgraph $M'$ of $M$ with $s(M')$ number of mixed cycles (length $\geq3$) gives $2^{s(M')}$ permutations $\pi$, since for each mixed cycle-component in $M'$ there are two ways of choosing the corresponding cycles in $\pi$. For a vertex $v_{k}\in V(M')$, we denote by $d_{k}=d(v_{k})=d_{M_{U}}(v_{k})$. Furthermore, if for some direction of a permutation $\pi$, a mixed cycle-component $C_{1}$ has value $\textbf{\textrm{i}}\prod_{v_{j}\in V(C_{1})}\frac{1}{d(v_{j})}$  (or $-\textbf{\textrm{i}}\prod_{v_{j}\in V(C_{1})}\frac{1}{d(v_{j})}$), then for the other direction $C_{1}$ has value $-\textbf{\textrm{i}}\prod_{v_{j}\in V(C_{1})}\frac{1}{d(v_{j})}$ (or $\textbf{\textrm{i}}\prod_{v_{j}\in V(C_{1})}\frac{1}{d(v_{j})}$) and vice verse. Thus, they cancel each other in the summation. In addition, if for some direction of a permutation $\pi$, $C_{1}$ has value $\prod_{v_{j}\in V(C_{1})}\frac{1}{d(v_{j})}$  (or $-\prod_{v_{j}\in V(C_{1})}\frac{1}{d(v_{j})}$), then for the other direction $C_{1}$ has the same value. For each edge-component $(kl)$ corresponds to the factors $(r_{h})_{kl}(r_{h})_{lk}$ has value $\frac{1}{\sqrt{d_{k}d_{l}}}\frac{1}{\sqrt{d_{l}d_{k}}}=\frac{1}{d_{k}d_{l}}$. For each arc-component $(kl)$ corresponds to the factors $(r_{h})_{kl}(r_{h})_{lk}$ has value $\frac{\textbf{\textrm{i}}\cdot(-\textbf{\textrm{i}})}{\sqrt{d_{k}d_{l}}\sqrt{d_{l}d_{k}}}=\frac{1}{d_{k}d_{l}}$.

Since $sgn(\pi)=(-1)^{n-c(M')}=(-1)^{r(M')}$ and each real spanning elementary subgraph $M'$ contributes $(-1)^{r(M')+l(M')}2^{s(M')}\prod_{v_{i}\in V(M')}\frac{1}{d_{M_{U}}(v_{i})}$ to the determinant of $R_{H}(M)$. This completes the proof.
\hfill$\blacksquare$

Now, we shall obtain a description of all the coefficients of the characteristic polynomial $P_{R_{H}}(M,x)$ of a mixed graph $M$.

\noindent\begin{theorem}\label{th:2.2}
For a mixed graph $M$, if the $R_{H}$-characteristic polynomial of $M$ is denoted by $P_{R_{H}}(M,x)=\mathrm{det}(xI-R_{H}(M))=x^{n}+a_{1}x^{n-1}+a_{2}x^{n-2}+\cdots+a_{n}$, then the coefficients of $P_{R_{H}}(M,x)$ are given by $$(-1)^{k}a_{k}=\sum_{M'}(-1)^{r(M')+l(M')}2^{s(M')}\prod_{v_{i}\in V(M')}\frac{1}{d_{M_{U}}(v_{i})},$$
where the summation is over all real  elementary subgraphs $M'$ with order $k$ of $M$, $r(M')=k-c(M')$, $c(M')$ denotes the number of components of $M'$, $l(M')$ denotes the number of negative mixed cycles of $M'$, $s(M')$ denotes the number of mixed cycles with length $\geq3$ in $M'$.
\end{theorem}

\noindent {\bf Proof.}
  The proof follows from Theorem \ref{th:2.1} and the fact that $(-1)^{k}a_{k}$ is the summation of determinants of all principal $k\times k$ submatrices of $R_{H}(M)$.  \hfill$\blacksquare$

 From Theorem \ref{th:2.2}, we can obtain the following results.

\noindent\begin{corollary}\label{co:2.3}
For a mixed graph $M$, let the $R_{H}$-characteristic polynomial of $M$ be denoted by $P_{R_{H}}(M,x)=|xI-R_{H}(M)|=x^{n}+a_{1}x^{n-1}+a_{2}x^{n-2}+\cdots+a_{n}$.
\begin{enumerate}[(1)]
  \item If $M$ is a mixed tree, then $$(-1)^{k}a_{k}=\sum_{M'}(-1)^{r(M')}\prod_{v_{i}\in V(M')}\frac{1}{d_{M_{U}}(v_{i})}.$$
  \item If $M$ is a mixed graph and its underlying graph $M_{U}$ is $r$ regular $(r\neq0)$, then $$(-1)^{k}a_{k}=\sum_{M'}(-1)^{r(M')+l(M')}2^{s(M')}\frac{1}{r^{k}}.$$
  \item If $M$ is a mixed bipartite graph, then all coefficients of $a_{odd}$ are equal to 0, and its spectrum is symmetry about 0.
\end{enumerate}
\end{corollary}

Note that if $M$ is a positive mixed graph, then for every real elementary subgraph $M'$ of $M$, we have
\begin{align*}
&(-1)^{r(M')+l(M')}2^{s(M')}\prod_{v_{i}\in V(M')}\frac{1}{d_{M_{U}}(v_{i})}\\
 &=(-1)^{r(M')}2^{s(M')}\prod_{v_{i}\in V(M')}\frac{1}{d_{M_{U}}(v_{i})}\\
 &=(-1)^{r(M_{U}')}2^{s(M_{U}')}\prod_{v_{i}\in V(M_{U}')}\frac{1}{d_{M_{U}}(v_{i})}.
\end{align*}

Then $P_{R_{H}}(M,x)=P_{R_{H}}(M_{U},x)$, that is to say

\noindent\begin{theorem}\label{th:2.4}
If $M$ is a positive mixed graph and $M_{U}$ be its underlying graph, then $Sp_{R_{H}}(M)=Sp_{R_{H}}(M_{U})$.
\end{theorem}

\section{Bounds on Hermitian-Randi\'{c} energy of mixed graphs}
In this section, we will give some bounds about Hermitian-Randi\'{c} energy of  mixed graphs.
First, we will give some properties of Hermitian-Randi\'{c} matrix of mixed graphs.
\noindent\begin{lemma}\label{le:3.1}
Let $M$ be a mixed graph of order $n\geq1$.
\begin{enumerate}[(1)]
  \item  $\mathcal{E}_{R_{H}}(M)=0$ if and only if $M\cong\overline{K}_{n}$.
  \item If $M=M_{1}\cup M_{2}\cup\cdots\cup M_{p}$, then $\mathcal{E}_{R_{H}}(M)=\mathcal{E}_{R_{H}}(M_{1})+\mathcal{E}_{R_{H}}(M_{2})+\cdots+\mathcal{E}_{R_{H}}(M_{p}).$
\end{enumerate}
\end{lemma}

From Lemma \ref{le:3.1}, we can obtain the following theorem.

\noindent\begin{theorem}\label{th:3.2}
Let $M$ be a mixed graph with vertex set $V(M)=\{v_{1},v_{2},\ldots,v_{n}\}$ and $d_{k}$ is the degree of $v_{k}$, $k=1,2,\ldots,n$. Let $H(M)$ and $R_{H}(M)$ be the Hermitian-adjacency matrix and Hermitian-Randi\'{c} matrix of $M$, respectively. If $M$ has isolated vertices, then $\emph{det}H(M)=\emph{det}R_{H}(M)=0$. If $M$ has no isolated vertices, then
$$\emph{det}R_{H}(M)=\frac{1}{d_{1}d_{2}\cdots d_{n}}\emph{det}H(M).$$
\end{theorem}

\noindent {\bf Proof.} If $M$ has $l$ isolated vertices, then $M=M'\cup \overline{K}_{l}$, where $M'$ has no isolated vertices. By Lemma \ref{le:3.1}, we have $Sp_{R_{H}}(M)=Sp_{R_{H}}(M')\cup\{0, l~times\}$ and an analogous relation holds for Hermitian-adjacency spectrum of $M$. That is, $H(M)$ and $R_{H}(M)$ have zero eigenvalues, therefore their determinants are equal to zero.
If $M$ has no isolated vertices, then $R_{H}(M)=D(M_{U})^{-\frac{1}{2}}H(M)D(M_{U})^{-\frac{1}{2}}$ is applicable,
where $D(M_{U})$ is the diagonal matrix of vertex degrees. The matrices $R_{H}(M)$ and $D(M_{U})^{-\frac{1}{2}}R_{H}(M)D(M_{U})^{\frac{1}{2}}$ are similar and thus have equal eigenvalues. We have
$$D(M_{U})^{-\frac{1}{2}}R_{H}(M)D(M_{U})^{\frac{1}{2}}=D(M_{U})^{-1}H(M),$$
and therefore, $$\textrm{det}R_{H}(M)=\textrm{det}[D(M_{U})^{-1}H(M)]=\textrm{det}D(M_{U})^{-1}\textrm{det}H(M).$$

So,
$$\textrm{det}R_{H}(M)=\frac{1}{d_{1}d_{2}\cdots d_{n}}\textrm{det}H(M).$$

This completes the proof.
\hfill$\blacksquare$

Similar to Theorem \ref{th:3.2}, we can obtain the following theorem.

\noindent\begin{theorem}\label{th:3.3}

If $M$ is a mixed graph with vertex set $V(M)=\{v_{1},v_{2},\ldots,v_{n}\}$ and its underlying graph $M_{U}$ is $r$ regular, then $\displaystyle\mathcal{E}_{R_{H}}(M)=\frac{1}{r}\mathcal{E}_{H}(M)$. In addition, if $r=0$, then $\mathcal{E}_{R_{H}}(M)=0$.
\end{theorem}

\noindent {\bf Proof.} If $r=0$, then $M$ is the graph that has no edges. Then all the entries of $R_{H}(M)$ are equal to 0, i.e., $R_{H}(M)=\emph{\textbf{0}}$. Similarly, $H(M)=\emph{\textbf{0}}$. Since all eigenvalues of the zero matrix are equal to 0. Hence, $\displaystyle\mathcal{E}_{R_{H}}(M)=\mathcal{E}_{H}(M)$=0.

If $r>0$, i.e., $M$ is regular of degree $r>0$. Then $d_{1}=d_{2}=\cdots=d_{n}=r$, where $d_{k}$ is the degree of $v_{k}$, $k=1,2,\ldots,n$. Hence
\begin{displaymath}
(r_{h})_{sk}=\left\{\
        \begin{array}{ll}
          \frac{1}{r},&v_{s}\leftrightarrow v_{k},\\
         \frac{\textbf{\textrm{i}}}{r},&v_{s}\rightarrow v_{k},\\
         \frac{-\textbf{\textrm{i}}}{r},&v_{k}\rightarrow v_{s},\\
         0,&\textrm{otherwise}.
        \end{array}
      \right.
\end{displaymath}

Which implies that $\displaystyle R_{H}(M)=\frac{1}{r}H(M)$.  Therefore, $\displaystyle\mu_{i}=\frac{1}{r}\lambda_{i}$, where $\mu_{i}$ is the eigenvalue of $R_{H}(M)$, and $\lambda_{i}$ is the eigenvalue of $H(M)$, for $i=1,2,\ldots,n$. Then this theorem follows from the definitions of $\mathcal{E}_{R_{H}}(M)$ and $\mathcal{E}_{H}(M)$. \hfill$\blacksquare$

Similar the results about the skew Randi\'{c} energy in \cite{GHL}, we can establish the following lower and upper bounds for the Hermitian-Randi\'{c} energy. First, we need the following theorem. Here and later by $\textrm{\textbf{I}}_{n}$ is denoted by the unit matrix of order $n$.

\noindent\begin{theorem}\label{th:3.4}
Let $M$ be a mixed graph of order $n$  and $\mu_{1}\geq\mu_{2}\geq\cdots\geq\mu_{n}$ be the Hermitian-Randi\'{c} spectrum of $R_{H}(M)$. Then $|\mu_{1}|=|\mu_{2}|=\cdots=|\mu_{n}|$ if and only if  there exists a constant $c=|\mu_{i}|^{2}$ for all $i$ such that $R^{2}_{H}(M)=c \emph{\textbf{I}}_{n}$.
\end{theorem}
\noindent {\bf Proof.} Let $\mu_{1}\geq\mu_{2}\geq\cdots\geq\mu_{n}$ be the Hermitian-Randi\'{c} spectrum of $R_{H}(M)$. Then there exists a unitary matrix $U$ such that
\begin{displaymath}
U^{\ast}R_{H}(M)U=U^{\ast}R_{H}(M)^{\ast}U=\left(
  \begin{array}{ccccccccccc}
        \mu_{1}         &          &\\
             &      \ddots         &\\
             &            &        \mu_{n}\\

  \end{array}
\right).
\end{displaymath}

So,
\begin{align*}
&|\mu_{1}|=|\mu_{2}|=\cdots=|\mu_{n}|\\
 \Leftrightarrow& ~U^{\ast}R_{H}(M)^{\ast}R_{H}(M)U=c \textrm{\textbf{I}}_{n}\\
 \Leftrightarrow& ~U(U^{\ast}R_{H}(M)^{\ast}R_{H}(M)U)U^{\ast}=cUU^{\ast}\\
 \Leftrightarrow& ~R_{H}(M)^{\ast}R_{H}(M)=c \textrm{\textbf{I}}_{n}\\
 \Leftrightarrow& ~R^{2}_{H}(M)=c \textrm{\textbf{I}}_{n},
\end{align*}
where $c$ is a constant and $c=|\mu_{i}|^{2}$ for all $i$.

This completes the proof.
\hfill$\blacksquare$

\noindent\begin{theorem}\label{th:3.5}
Let $M$ be a mixed graph of order $n$ and $\mu_{1}\geq\mu_{2}\geq\cdots\geq\mu_{n}$ be the Hermitian-Randi\'{c} spectrum of $R_{H}(M)$. Let $M_{U}$ be the underlying graph of $M$, $p=|\emph{det}R_{H}(M)|$. Then $$\sqrt{2R_{-1}(M_{U})+n(n-1)p^{\frac{2}{n}}}\leq\mathcal{E}_{R_{H}}(M)\leq\sqrt{2nR_{-1}(M_{U})},$$
with equalities hold both in the lower bound and upper bound if and only if there exists a constant $c=|\mu_{i}|^{2}$ for all $i$ such that $R^{2}_{H}(M)=c \emph{\textbf{I}}_{n}$.
\end{theorem}

\noindent {\bf Proof.} Let $\{\mu_{1},\mu_{2},\ldots,\mu_{n}\}$ be the Hermitian-Randi\'{c} spectrum of $M$, where $\mu_{1}\geq\mu_{2}\geq\cdots\geq\mu_{n}$. Since
 $\sum_{j=1}^{n}\mu_{j}^{2}=tr(R_{H}^{2}(M))=\sum_{j=1}^{n}\sum_{k=1}^{n}(r_{h})_{jk}(r_{h})_{kj}=\sum_{j=1}^{n}\sum_{k=1}^{n}(r_{h})_{jk}\overline{(r_{h})}_{jk}$
$=\sum_{j=1}^{n}\sum_{k=1}^{n}|(r_{h})_{jk}|^{2}=2R_{-1}(M_{U}),$
where $R_{-1}(M_{U})=\sum_{v_{j}v_{k}\in E(M_{U})}\frac{1}{d_{j}d_{k}}$ (unordered).

Applying the Cauchy-Schwartz inequality we have
$$\mathcal{E}_{R_{H}}(M)=\sum^{n}_{j=1}|\mu_{j}|\leq\sqrt{\sum^{n}_{j=1}|\mu_{j}|^{2}}\cdot\sqrt{n}=\sqrt{2nR_{-1}(M_{U})}.$$

On the other hand,
$$|\mathcal{E}_{R_{H}}(M)|^{2}=\left(\sum^{n}_{j=1}|\mu_{j}|\right)^{2}=\sum^{n}_{j=1}|\mu_{j}|^{2}+\sum_{1\leq i\neq j\leq n}|\mu_{i}||\mu_{j}|.$$

By using arithmetic geometric average inequality, we can get that
$$|\mathcal{E}_{R_{H}}(M)|^{2}=\sum^{n}_{j=1}|\mu_{j}|^{2}+\sum_{1\leq i\neq j\leq n}|\mu_{i}||\mu_{j}|\geq2R_{-1}(M_{U})+n(n-1)p^{\frac{2}{n}}.$$

Therefore, we can obtain the lower bound on Hermitian-Randi\'{c} energy,
 $$\mathcal{E}_{R_{H}}(M)\geq\sqrt{2R_{-1}(M_{U})+n(n-1)p^{\frac{2}{n}}}.$$

From the Cauchy-Schwartz inequality and arithmetic geometric average inequality, we know that the equalities hold both in the lower bound and upper bound if and only if $|\mu_{1}|=|\mu_{2}|=\cdots=|\mu_{n}|$, i.e., there exists a constant $c=|\mu_{i}|^{2}$ for all $i$ such that $R^{2}_{H}(M)=c \textrm{\textbf{I}}_{n}$.

 This completes the proof.
\hfill$\blacksquare$

\noindent\begin{corollary}\label{co:3.6}
Let $M$ be a mixed graph and its underlying graph $M_{U}$ be $r$ $(\neq0)$ regular and $E(M_{U})=m$. Let $\mu_{1}\geq\mu_{2}\geq\cdots\geq\mu_{n}$ be the Hermitian-Randi\'{c} spectrum of $R_{H}(M)$. Then
$$\sqrt{\frac{n}{r}+n(n-1)p^{\frac{2}{n}}}\leq\mathcal{E}_{R_{H}}(M)\leq\frac{n\sqrt{r}}{r},$$ where $p=|\emph{det}R_{H}(M)|,$
with equalities hold both in the lower bound and upper bound if and only if $\displaystyle\frac{1}{r}=|\mu_{i}|^{2}$ for all $i$ such that $\displaystyle R^{2}_{H}(M)=\frac{1}{r} \emph{\textbf{I}}_{n}$.
\end{corollary}
\noindent {\bf Proof.} If $M$ is a mixed graph and its underlying graph $M_{U}$ is $r$ regular, then $\displaystyle R_{-1}(M_{U})=\frac{m}{r^{2}}$ and $2m=nr$. By Theorems \ref{th:3.4} and \ref{th:3.5}, we can obtain the results.
\hfill$\blacksquare$

\noindent\begin{lemma}\label{le:3.7}\cite{LXYYS}
Let $G$ be a graph of order $n$ with no isolated vertices. Then
$$\frac{n}{2(n-1)}\leq R_{-1}(G)\leq\lfloor\frac{n}{2}\rfloor,$$
with equality in the lower bound if and only if $G$ is a complete graph, and equality in the
upper bound if and only if either
\begin{enumerate}[(1)]
  \item $n$ is even and $G$ is the disjoint union of $n/2$ paths of length 1, or
  \item $n$ is odd and $G$ is the disjoint union of $(n-3)/2$ paths of length 1 and one path of length 2.
\end{enumerate}
\end{lemma}

Combining with Theorem \ref{th:3.5} and Lemma \ref{le:3.7}, we can get upper and lower bounds for the Hermitian-Randi\'{c} energy by replacing $R_{-1}(M_{U})$ with other parameters. We now give bounds of Hermitian-Randi\'{c} energy of a mixed graph with respect to its order.

\noindent\begin{theorem}\label{th:3.8}
Let $M$ be a mixed graph of order $n\geq3$ without isolated vertices and $M_{U}$ be its underlying graph. Let $\mu_{1}\geq\mu_{2}\geq\cdots\geq\mu_{n}$ be the Hermitian-Randi\'{c} spectrum of $R_{H}(M)$. Then
$$\sqrt{\frac{2n}{n-1}}\leq \mathcal{E}_{R_{H}}(M)\leq n.$$
The equality in the upper bound holds if and only if $n$ is even and $M_{U}$ is the disjoint union of $n/2$ paths of length 1. The equality in the lower bound holds if and only if $M_{U}$ is a complete graph and $\mu_{1}=-\mu_{n}\neq0$, $\mu_{j}=0$, $j=2,\ldots,n-1$.
\end{theorem}

\noindent {\bf Proof.} Let $R_{H}(M)$ be the Hermitian-Randi\'{c} matrix of $M$ and $\mu_{1}\geq\mu_{2}\geq\cdots\geq\mu_{n}$ be the Hermitian-Randi\'{c} spectrum of $R_{H}(M)$.

For the upper bound, combining with Lemma \ref{le:3.7} and $\mathcal{E}_{R_{H}}(M)\leq\sqrt{2nR_{-1}(M_{U})}$ of Theorem \ref{th:3.5}, we have
$$\mathcal{E}_{R_{H}}(M)\leq\sqrt{2nR_{-1}(M_{U})}\leq\sqrt{2n\lfloor\frac{n}{2}\rfloor}\leq n.$$
From Theorem \ref{th:3.5} and Lemma \ref{le:3.7}, we know that the equality in the upper bound holds if and
only if $n$ is even, $M_{U}$ is the graph described in Lemma \ref{le:3.7}~(1), and
$|\mu_{1}|=|\mu_{2}|=\cdots=|\mu_{n}|$, that is, we can obtain the upper bound when $n$ is even and $M_{U}$
is the disjoint union of $n/2$ paths of length 1.

For the lower bound, since the sum of the diagonal entries of $R_{H}(M)$ is 0, i.e., $\sum^{n}_{k=1}\mu_{k}=0$. Then
\begin{align*}
\left(\sum^{n}_{k=1}\mu_{k}\right)\left(\sum^{n}_{l=1}\mu_{l}\right)&=\sum^{n}_{k=1}\mu_{k}^{2}+\sum_{1\leq k\neq l\leq n}\mu_{k}\mu_{l}\\
 &=\sum^{n}_{k=1}\mu_{k}^{2}+2\sum_{k<l}\mu_{k}\mu_{l}\\
 &=2R_{-1}(M_{U})+2\sum_{k<l}\mu_{k}\mu_{l}\\
 &=0.
\end{align*}

Hence, $\displaystyle\sum_{k<l}\mu_{k}\mu_{l}=-R_{-1}(M_{U})$.

From the definition of the Hermitian-Randi\'{c} energy of a mixed graph, we have
\begin{align*}
\mathcal{E}^{2}_{R_{H}}(M)&=\left(\sum^{n}_{k=1}|\mu_{k}|\right)^{2}\\
 &=\sum^{n}_{k=1}\mu_{k}^{2}+\sum_{1\leq k\neq l\leq n}|\mu_{k}\mu_{l}|\\
 &=2R_{-1}(M_{U})+2\sum_{k<l}|\mu_{k}\mu_{l}|\\
 &\geq2R_{-1}(M_{U})+2|\sum_{k<l}\mu_{k}\mu_{l}|\\
 &=4R_{-1}(M_{U}).
\end{align*}

Combining with Lemma \ref{le:3.7}, we have
$$\displaystyle\mathcal{E}_{R_{H}}(M)\geq2\sqrt{R_{-1}(M_{U})}\geq2\sqrt{\frac{n}{2(n-1)}}=\sqrt{\frac{2n}{n-1}}.$$

So, $$\displaystyle\sqrt{\frac{2n}{n-1}}\leq\mathcal{E}_{R_{H}}(M)\leq n.$$

From the proof above and Lemma \ref{le:3.7}, we know that the equality in the lower bound holds if and only if
$M_{U}$ is a complete graph and $\mu_{k}\mu_{l}\geq 0$ or $\mu_{k}\mu_{l}\leq 0$,  for all $1\leq k< l\leq n$.  Note that $\sum_{k=1}^{n}\mu_{k}=0$ and $M$ has no isolated vertices, so the former case can not happen. Hence, the equality in the lower bound holds if and only if $M_{U}$ is a complete graph and $\mu_{1}=-\mu_{n}\neq0$, $\mu_{j}=0$, $j=2,\ldots,n-1$.

This completes the proof.
\hfill$\blacksquare$

\noindent\begin{remark}\label{re:3.9}
  \rm It should be pointed out that when $M$ is a complete mixed graph, its Hermitian-Randi\'{c} spectrum is not unique. For example, let $M_{U}=K_{3}$, if all edges of $E(M)$ are oriented, then we have $\mu_{1}=-\mu_{3}=\frac{\sqrt{3}}{2}$, $\mu_{2}=0$, then we can obtain the lower bound in Theorem \ref{th:3.8}. If there exist some edges of  $E(M)$ are undirected, then we can not obtain the lower bound in Theorem \ref{th:3.8}. For example, if $(r_{h})_{12}=(r_{h})_{32}=\frac{i}{2}$, $(r_{h})_{13}=\frac{1}{2}$, then $\mu_{1}=1$ and $\mu_{2}=\mu_{3}=-\frac{1}{2}$.  Hence, the problem of determining all complete mixed graphs for which the lower bound in Theorem \ref{th:3.8} is attained appears
to be somewhat more difficult.
\end{remark}

To deduce more bounds on $\mathcal{E}_{R_{H}}(M)$, the following lemma is needed.

\noindent\begin{lemma}\label{le:3.10}\cite{SD}
Let $\textbf{x},~\textbf{y}\in \mathbb{R}^{n}$ and let $A(\textbf{x})=\frac{1}{n}\sum^{n}_{i=1}x_{i}$, $A(\textbf{y})=\frac{1}{n}\sum^{n}_{i=1}y_{i}$. If $\phi\leq x_{i}\leq\Phi$ and $\gamma\leq y_{i}\leq\Gamma$, then
$$\left|\frac{1}{n}\sum^{n}_{i=1}x_{i}y_{i}-\frac{1}{n^{2}}\sum^{n}_{i=1}x_{i}\sum^{n}_{i=1}y_{i}\right|\leq\sqrt{(\Phi-A(x))(A(x)-\phi)(\Gamma-A(y))(A(y)-\gamma)}.$$
\end{lemma}

Now we turn to new bounds on $\mathcal{E}_{R_{H}}(M)$.

\noindent\begin{theorem}\label{th:3.11}
Let $M$ be a mixed graph of order $n$ and $M_{U}$ be its underlying graph.  Let $\mu_{1}\geq\mu_{2}\geq\cdots\geq\mu_{n}$ be the Hermitian-Randi\'{c} spectrum of $R_{H}(M)$. Then
\begin{equation}\label{1}
\mathcal{E}_{R_{H}}(M)\geq\frac{2R_{-1}(M_{U})+n\alpha\beta}{\alpha+\beta},
\end{equation}
where $\alpha=\min \limits_{1\leq i\leq n}\{|\mu_{i}|\}$, $\beta=\max\{\mu_{1},|\mu_{n}|\}$.

\end{theorem}

\noindent {\bf Proof.} Note that
\begin{equation}\label{2}
\begin{array}{ll}
  \mathcal{E}^{2}_{R_{H}}(M)=\left(\sum^{n}_{j=1}|\mu_{j}|\right)^{2}
  =\sum^{n}_{j=1}|\mu_{j}|^{2}+\sum_{1\leq i\neq j\leq n}|\mu_{i}||\mu_{j}|\\
~~~~~~~~~~~~~~=2R_{-1}(M_{U})+\sum_{1\leq i\neq j\leq n}|\mu_{i}||\mu_{j}|.
\end{array}
\end{equation}

Let $S=\sum_{1\leq i\neq j\leq n}|\mu_{i}||\mu_{j}|$, $x_{i}=|\mu_{i}|$ and $y_{i}=\mathcal{E}_{R_{H}}(M)-|\mu_{i}|$, $i=1,2,\ldots,n$.
Then $S=\sum^{n}_{i=1}x_{i}y_{i}$.

From the definitions of $\alpha$ and $\beta$, we have $\alpha\leq x_{i}\leq\beta$ and $\mathcal{E}_{R_{H}}(M)-\beta\leq y_{i}\leq\mathcal{E}_{R_{H}}(M)-\alpha$. In addition, let $A(\textbf{x})=\frac{1}{n}\sum^{n}_{i=1}x_{i}=\frac{\mathcal{E}_{R_{H}}(M)}{n}$ and $ A(\textbf{y})=\frac{1}{n}\sum^{n}_{i=1}y_{i}=\frac{(n-1)\mathcal{E}_{R_{H}}(M)}{n}$. Hence by Lemma \ref{le:3.10}, we have
\begin{align*}
\left|\frac{S}{n}-\frac{(n-1)\mathcal{E}^{2}_{R_{H}}(M)}{n^{2}}\right| &\leq\sqrt{\left(\beta-\frac{\mathcal{E}_{R_{H}}(M)}{n}\right)\left(\frac{\mathcal{E}_{R_{H}}(M)}{n}-\alpha\right)}\\
\\
&\cdot\sqrt{\left[\mathcal{E}_{R_{H}}(M)-\alpha-\frac{(n-1)\mathcal{E}_{R_{H}}(M)}{n}\right]}\\
&\cdot\sqrt{\left[\frac{(n-1)\mathcal{E}_{R_{H}}(M)}{n}-(\mathcal{E}_{R_{H}}(M)-\beta)\right]}\\
 &=\sqrt{\left(\beta-\frac{\mathcal{E}_{R_{H}}(M)}{n}\right)^{2}\left(\frac{\mathcal{E}_{R_{H}}(M)}{n}-\alpha\right)^{2}}.
\end{align*}

It follows that
$$S\geq\mathcal{E}^{2}_{R_{H}}(M)+n\alpha\beta-(\alpha+\beta)\mathcal{E}_{R_{H}}(M).$$ This together with (\ref{2}) implies that
$$\mathcal{E}^{2}_{R_{H}}(M)=2R_{-1}(M_{U})+S\geq2R_{-1}(M_{U})+\mathcal{E}^{2}_{R_{H}}(M)+n\alpha\beta-(\alpha+\beta)\mathcal{E}_{R_{H}}(M).$$

So,
$$\mathcal{E}_{R_{H}}(M)\geq\frac{2R_{-1}(M_{U})+n\alpha\beta}{\alpha+\beta}.$$

This completes the proof.
\hfill$\blacksquare$

Note that the right-hand side of (\ref{1}) is a  non-decreasing function on $\alpha\geq0$. Combining with Theorem \ref{th:3.4}, we have the following corollary.

\noindent\begin{corollary}\label{co:3.12}
Let $M$ be a mixed graph of order $n$ and $M_{U}$ be its underlying graph. Let $\mu_{1}\geq\mu_{2}\geq\cdots\geq\mu_{n}$ be the Hermitian-Randi\'{c} spectrum of $R_{H}(M)$. Then
$$\mathcal{E}_{R_{H}}(M)\geq\frac{2R_{-1}(M_{U})}{\beta},$$
where $\beta=\max\{\mu_{1},|\mu_{n}|\}$. The equality holds if and only if $R^{2}_{H}(M)=c \emph{\textbf{I}}_{n}$, where $c$ is a constant such that $|\mu_{i}|^{2}=c$ for all $i$.

\end{corollary}

In particular, if $M$ is a connected mixed bipartite graph, then we have the following theorem.

\noindent\begin{theorem}\label{th:3.13}
Let $M$ be a connected mixed bipartite graph of order $n$  and $M_{U}$ be its underlying graph. Let $\mu_{1}\geq\mu_{2}\geq\cdots\geq\mu_{n}$ be the Hermitian-Randi\'{c} spectrum of $R_{H}(M)$. Then
\begin{equation}\label{3}
\mathcal{E}_{R_{H}}(M)\geq2\left(\frac{R_{-1}(M_{U})+\lfloor\frac{n}{2}\rfloor\alpha\mu_{1}}{\alpha+\mu_{1}}\right),
\end{equation}
where $\alpha=\min \limits_{1\leq i\leq \lfloor\frac{n}{2}\rfloor}\{|\mu_{i}|\}$.

\end{theorem}

\noindent {\bf Proof.} Note that $M_{U}$ is a bipartite graph. By Corollary \ref{co:2.3} (3), we have $\mu_{i}=-\mu_{n+1-i}$ and $\mu_{i}\geq0$ for $i=1,2,\ldots,\lfloor\frac{n}{2}\rfloor$. Therefore,
\begin{equation}
\begin{array}{ll}\label{4}
 \mathcal{E}^{2}_{R_{H}}(M)=\left(2\sum^{\lfloor\frac{n}{2}\rfloor}_{i=1}\mu_{i}\right)^{2}
  =4(\sum^{\lfloor\frac{n}{2}\rfloor}_{i=1}\mu_{i}^{2}+\sum_{1\leq i\neq j\leq \lfloor\frac{n}{2}\rfloor}\mu_{i}\mu_{j})\\
~~~~~~~~~~~~~~=4R_{-1}(M_{U})+4\sum_{1\leq i\neq j\leq \lfloor\frac{n}{2}\rfloor}\mu_{i}\mu_{j}.
\end{array}
\end{equation}

Let $T=\sum_{1\leq i\neq j\leq \lfloor\frac{n}{2}\rfloor}\mu_{i}\mu_{j}$, $x_{i}=\mu_{i}$ and $y_{i}=\frac{\mathcal{E}_{R_{H}}(M)}{2}-\mu_{i}$, $i=1,2,\ldots,\lfloor\frac{n}{2}\rfloor$.
Then $T=\sum^{\lfloor\frac{n}{2}\rfloor}_{i=1}x_{i}y_{i}$.

From the definition of $\alpha$, we have $\alpha\leq x_{i}\leq\mu_{1}$ and $\frac{\mathcal{E}_{R_{H}}(M)}{2}-\mu_{1}\leq y_{i}\leq\frac{\mathcal{E}_{R_{H}}(M)}{2}-\alpha$. In addition, let $ A(\textbf{x})=\frac{1}{\lfloor\frac{n}{2}\rfloor}\sum^{\lfloor\frac{n}{2}\rfloor}_{i=1}x_{i}=\frac{\mathcal{E}_{R_{H}}(M)}{2\lfloor\frac{n}{2}\rfloor}$ and $ A(\textbf{y})=\frac{1}{\lfloor\frac{n}{2}\rfloor}\sum^{\lfloor\frac{n}{2}\rfloor}_{i=1}y_{i}=\frac{(\lfloor\frac{n}{2}\rfloor-1)\mathcal{E}_{HR}(M)}{2\lfloor\frac{n}{2}\rfloor}$. Hence by Lemma \ref{le:3.10}, we have
\begin{align*}
\left|\frac{T}{\lfloor\frac{n}{2}\rfloor}-\frac{(\lfloor\frac{n}{2}\rfloor-1)\mathcal{E}^{2}_{R_{H}}(M)}{4\lfloor\frac{n}{2}\rfloor^{2}}\right| &\leq\sqrt{\left(\mu_{1}-\frac{\mathcal{E}_{R_{H}}(M)}{2\lfloor\frac{n}{2}\rfloor}\right)^{2}\left(\frac{\mathcal{E}_{R_{H}}(M)}{2\lfloor\frac{n}{2}\rfloor}-\alpha\right)^{2}}.
\end{align*}

It follows that
$$T\geq\frac{\mathcal{E}^{2}_{R_{H}}(M)}{4}+\lfloor\frac{n}{2}\rfloor\alpha\mu_{1}-(\alpha+\mu_{1})\frac{\mathcal{E}_{R_{H}}(M)}{2}.$$ This together with (\ref{4}) implies that
$$\mathcal{E}^{2}_{R_{H}}(M)=4R_{-1}(M_{U})+4T\geq4R_{-1}(M_{U})+\mathcal{E}^{2}_{R_{H}}(M)+4\lfloor\frac{n}{2}\rfloor\alpha\mu_{1}-2(\alpha+\mu_{1})\mathcal{E}_{R_{H}}(M).$$

So,
$$\mathcal{E}_{R_{H}}(M)\geq2\left(\frac{R_{-1}(M_{U})+\lfloor\frac{n}{2}\rfloor\alpha\mu_{1}}{\alpha+\mu_{1}}\right).$$

This completes the proof.
\hfill$\blacksquare$

Note that the right-hand side of (\ref{3}) is a  non-decreasing function on $\alpha\geq0$. Combining with Theorem \ref{th:3.4}, we have the following corollary.

\noindent\begin{corollary}\label{co:3.14}
Let $M$ be a connected mixed bipartite graph of order $n$ and $M_{U}$ be its underlying graph. Let $\mu_{1}\geq\mu_{2}\geq\cdots\geq\mu_{n}$ be the Hermitian-Randi\'{c} spectrum of $R_{H}(M)$. Then
$$\mathcal{E}_{R_{H}}(M)\geq\frac{2R_{-1}(M_{U})}{\mu_{1}},$$
the equality holds if and only if  $R^{2}_{H}(M)=c \emph{\textbf{I}}_{n}$, where $c$ is a constant such taht $|\mu_{i}|^{2}=c$ for all $i$.

\end{corollary}

\section{Hermitian-Randi\'{c} energy of trees}
 In \cite{CRST}, the authors proved that the skew energy of a directed tree is independent of its orientation. In \cite{GHL}, the authors showed that the skew Randi\'{c} energy of a directed tree has the same result. In this section, we will show that the Hermitian-Randi\'{c} energy also has the same result.
In the beginning of this section, we first to characterize the mixed graphs with cut-edge.
\noindent\begin{theorem}\label{th:4.1}
Let $M$ be a mixed graph of order $n$ and $e=uv$ is an edge of $M$. If $uv$ is a cut-edge of $M_{U}$, where $M_{U}$ is the underlying graph of $M$. Then the spectrum and energy of $R_{H}(M)$ are unchanged when the edge $uv$ is replaced with a single arc $uv$ or $vu$ and vice versa.
\end{theorem}
\noindent {\bf Proof.} Let $e=uv$ be a cut-edge of $M_{U}$. Suppose $M_{1}$ is the graph obtained from $M$ by replacing the edge $uv$ with the arc $uv$ or $vu$. Let $M'$ and $M_{1}'$ be real elementary subgraph of order $k$ of $M$ and $M_{1}$, respectively. If  $M'$  do not contain the cut-edge $uv$, then $M'$ is also a real elementary subgraph of $M_{1}$, that is, $M'=M_{1}'$. So we have
\begin{equation}\label{5}
(-1)^{r(M')+l(M')}2^{s(M')}\prod_{v_{i}\in V(M')}\frac{1}{d_{M_{U}}(v_{i})}=(-1)^{r(M_{1}')+l(M_{1}')}2^{s(M_{1}')}\prod_{v_{i}\in V(M_{1}')}\frac{1}{d_{M_{U}}(v_{i})}.
\end{equation}

 If $M'$ contains the cut-edge $uv$, then there is a real elementary subgraph $M_{1}'$ of $M_{1}$ only differ from $M'$ on $uv$. Since $uv$ is a cut-edge of $M_{U}$, $uv$ is not contained in any cycles of $M_{U}$. Hence we have
 \begin{equation}\label{6}
 (-1)^{r(M')+l(M')}2^{s(M')}\prod_{v_{i}\in V(M')}\frac{1}{d_{M_{U}}(v_{i})}=(-1)^{r(M_{1}')+l(M_{1}')}2^{s(M_{1}')}\prod_{v_{i}\in V(M_{1}')}\frac{1}{d_{M_{U}}(v_{i})}.
 \end{equation}

Combining with (\ref{5}) and (\ref{6}), we have
\begin{align*}
  &a_{k}(M)-a_{k}(M_{1}) \\
  &=\sum_{M'}(-1)^{r(M')+l(M')}2^{s(M')}\prod_{v_{i}\in V(M')}\frac{1}{d_{M_{U}}(v_{i})}\\
  &-\sum_{M_{1}'}(-1)^{r(M_{1}')+l(M_{1}')}2^{s(M_{1}')}\prod_{v_{i}\in V(M_{1}')}\frac{1}{d_{M_{U}}(v_{i})}\\
  &=0,
\end{align*}
for any integer $k$.
Thus $Sp_{R_{H}}(M)=Sp_{R_{H}}(M_{1})$. Moreover, $\mathcal{E}_{{R_{H}}}(M)=\mathcal{E}_{{R_{H}}}(M_{1})$.

Similarly, we can prove that $Sp_{R_{H}}(M)=Sp_{R_{H}}(M_{2})$ and $\mathcal{E}_{{R_{H}}}(M)=\mathcal{E}_{{R_{H}}}(M_{2})$,
where $M_{2}$ is the mixed graph obtained from $M$ by replacing the arc $uv$ or $vu$ with the edge $uv$.
\hfill$\blacksquare$

Thus, the Hermitian-Randi\'{c} spectrum and Hermitian-Randi\'{c} energy are invariants when reversing the cut-arc's orientation or unorienting it or orienting an undirected cut-edge. By applying Theorem \ref{th:4.1},
we can obtain the following corollaries.

\noindent\begin{corollary}\label{co:4.2}
Let $T$ be a mixed tree of order $n$, and $T'$ be the mixed tree obtained from $T$ by reversing the orientations of
all the arcs incident with a particular vertex of $T$. Then
$\mathcal{E}_{R_{H}}(T)=\mathcal{E}_{R_{H}}(T')$.
\end{corollary}

\noindent\begin{corollary}\label{co:4.3}
Let $T$ be a mixed tree and $T_{U}$ be its underlying graph. Then
\begin{enumerate}[(1)]
  \item The Hermitian-Randi\'{c} energy of $T$ is independent of its orientation of the arc set.
  \item The Hermitian-Randi\'{c} energy of $T$ is the same as the Randi\'{c} energy of $T_{U}$.
\end{enumerate}
\end{corollary}

\end{document}